# Sharpening *Primes is in P* for a large family of numbers.

Pedro Berrizbeitia*

August 20, 2018

## Abstract

We present algorithms that are deterministic primality tests for a large family of integers, namely, integers $n \equiv 1 \pmod{4}$ for which an integer $a$ is given such that the Jacobi symbol $(\frac{a}{n}) = -1$, and integers $n \equiv -1 \pmod{4}$ for which an integer $a$ is given such that $(\frac{a}{n}) = (\frac{1-a}{n}) = -1$. The algorithms we present run in $2^{-\min(k,[2\log\log n])}\tilde{O}(\log n)^6$ time, where $k = \nu_2(n-1)$ is the exact power of 2 dividing $n-1$ when $n \equiv 1 \pmod{4}$ and $k = \nu_2(n+1)$ if $n \equiv -1 \pmod{4}$. The complexity of our algorithms improves up to $\tilde{O}(\log n)^4$ when $k \geq [2\log\log n]$. We also give tests for more general family of numbers and study their complexity.

## 1 Introduction

On August, 2002, Manindra Agrawal, Neeraj Kayal and Nitin Saxena published an important paper titled *Primes is in P* [3]. They produced an algorithm, now called the AKS algorithm, that determines whether a given number $n$ is prime or composite and that runs in polynomial time. This remarkable achievement gives a positive answer to the most important question in the general theory of primality testing. In fact, they prove that the AKS algorithm runs in $\tilde{O}((\log n)^{12})$ time, where $\tilde{O}(f(x))$ denotes $O(f(x)poly(\log f(x)))$.

In this paper we present algorithms that run faster than the AKS algorithm and are deterministic primality tests for a large family of integers, namely integers $n \equiv 1 \pmod{4}$ for which an integer $a$ is given such that the Jacobi symbol $(\frac{a}{n}) = -1$, and integers $n \equiv -1 \pmod{4}$ for which an integer $a$ is given such that $(\frac{a}{n}) = (\frac{1-a}{n}) = -1$. The algorithms we present run in $2^{-\min(k,[2\log\log n])}\tilde{O}(\log n)^6$ time, where $k = \nu_2(n-1)$ is the exact power of 2 dividing $n-1$ when $n \equiv 1 \pmod{4}$, and $k = \nu_2(n+1)$ if $n \equiv -1 \pmod{4}$. In particular, the running time of our algorithms improves up to $\tilde{O}(\log n)^4$ if the value of $k \geq [2\log\log n]$. If $n$ is a large enough prime, then we show that our algorithm for the case $n \equiv 1 \pmod{4}$ runs, in the worst case when $k = 2$, at least $2^{11}$ times faster than the best possible running time for the AKS algorithm. This advantage in running time increases with the value of $k$. For the case $n \equiv -1 \pmod{4}$ we get the same result using $2^9$ instead of $2^{11}$.

*Departamento de Matemáticas Puras y Aplicadas. Universidad Simón Bolívar, Apdo. 89000, Caracas 1080-A, Venezuela. email: *pedrob@usb.ve*

The first major breakthrough in the general theory of Primality Testing was achieved by Adleman, Pomerance and Rumely in 1983 [4], who gave a deterministic primality test running in $(\log n)^{O(\log \log \log n)}$ time. This algorithm was later improved and implemented by Cohen and Lenstra [8]. It is known in the literature as the APRCL algorithm. In [3] the authors present a brief summary of the main contributions to this general theory prior to AKS. They describe the contributions of Goldwasser and Kilian [10], of Atkin [1], and of Adleman an Huang [2].

The theory of primality testing for restricted families of numbers had an earlier start. The first and most famous "modern" algorithm is the Lucas-Lehmer Test [12]. It is an algorithm that runs in $\tilde{O}((\log n)^2)$ time to determine whether a Mersenne number (a number of the form $2^p - 1$, $p$ prime) is prime or composite. Proth [13] enlarged the family of numbers for which a primality test that runs in $\tilde{O}((\log n)^2)$ exists. The Proth Test applies to all numbers $n$ such that $\nu_2(n-1) > (1/2)\log n$ (by $\log n$ we always mean log to the base 2) provided an integer $a$ is given for which the Jacobi Symbol $(\frac{a}{n}) = -1$. Usually such an integer $a$ can easily be found using the quadratic reciprocity law; thus, the Proth test becomes deterministic for a large proportion of, though not all the numbers $n$ satisfying $\nu_2(n-1) > (1/2)\log n$. Later, the Lucas-Lehmer Test was also extended to all numbers $n \equiv -1 \pmod{4}$, such that $\nu_2(n+1) > (1/2)\log n$ for which an integer $a$ is given such that $(\frac{a}{n}) = (\frac{1-a}{n}) = -1$. In a series of papers starting around 1970, Hugh Williams and collaborators extended these tests to numbers satisfying $\nu_p(n \pm 1) > (1/2)\log n$, where $p$ is a prime, provided there is a prime $q$, $q \equiv 1 \pmod{p}$, such that $n$ is not a $p$-th power modulo $q$, and gave many concrete implementations and tables of primes. Further extensions of Williams results can be found in [7]. The book of Williams [15] is a good source for studying many of these results and the history of this subject.

Our paper links the two approaches described above: the general approach and the one for restricted families of numbers. We still need an integer $a$ satisfying the Jacobi Symbol condition, but we no longer impose any condition on $\nu_2(n-1)$ or on $\nu_2(n+1)$. Thus the tests can be implemented for a set of numbers of density arbitrarily near 1. The link is also evidenced by the fact that the complexity of the tests we give decreases as the value of $k$ increases.

As mentioned before, the authors of [3] were able to prove that the AKS algorithm runs in $\tilde{O}((\log n)^{12})$ time, but they believe (and have strong evidence to support this belief), that it actually runs in $\tilde{O}((\log n)^6)$ time. In fact they prove that this would be the case if a widely believed conjecture on the density of the Sophie-Germain primes is true. The main step of their algorithm (the step that determines the complexity) consists in verifying that
$$(m+x)^n \equiv m + x^n \pmod{n,\ x^r - 1}$$

for $m = 1$ to $2\sqrt{r}\log n$ where $r$ is a prime with specific properties ($r - 1$ has a prime divisor $q > 4\sqrt{r}\log n$ which divides the order of $n$ modulo $r$). They prove such prime r exists in the interval $(64(\log n)^2, c(\log n)^6)$ for some constant $c$. They are able to prove that $r < c(\log n)^6$ by making clever use of a result in analytic number theory on the density of primes. But they believe that such $r$ is actually of size $O(\log n)^2$ (and prove this under the assumption of the Sophie Germain prime density conjecture mentioned above). The lower bound on $r$ implies that the AKS algorithm runs in at least $\tilde{O}((\log n)^6)$ when n is a prime. The upper bound implies that it runs in at most $\tilde{O}((\log n)^{12})$. According to Bernstein [5] Lenstra was able to prove that in fact such $r$ is $O(\log n)^4$, hence proved that the complexity of AKS is at most $\tilde{O}((\log n)^8)$. He also showed that $r$ need not be a prime, but that could be any number such that $n$ is a primitive root modulo $r$.

In the case $n \equiv 1 \pmod 4$, and assuming an integer $a$ is given such that $\left(\frac{a}{n}\right) = -1$ the two key observations in our paper are:

1. It is enough to verify
$$(1 + mx)^n \equiv 1 + mx^n \pmod{n, \ x^{2^s} - a}$$
where $s = [2\log\log n]$ (hence $2^s < (\log n)^2$). Since $2^s$ is smaller than $r$ (in fact is it at least 64 times smaller than $r$) then each of these verifications for different values of $m$ are faster than the verification of the analogous step in the AKS algorithm.

2. These verifications only have to be done for $2^{\max(s-k,0)}$ different values of $m$, where $k = \nu_2(n-1)$. We will see this in detail within the proof of Theorem 3.1 and 4.1, but we point out here the crucial fact, namely, that some of the conjugates of the monomial $1 + mx^n$ in the corresponding finite field are also monomials satisfying the same congruence. So, each iteration of our test produces $2^{\min(s,k)}$ different monomials satisfying the congruence.

These two facts together allow us to give a more efficient primality test for those numbers and such that its efficiency improves with the value of $k$ up to a certain limit ($[2\log\log n]$). For numbers $n \equiv -1 \pmod 4$ we were able to obtain similar results.

In Section 2 we define the notation and give some elementary but necessary results on the theory of finite fields. In Section 3 we present the algorithm for the case $n \equiv 1 \pmod 4$, we prove the validity of the algorithm and study its complexity. In Section 4 we do the same for the case $n \equiv -1 \pmod 4$. In this case our algorithm runs around 4 times slower than the one given in the previous section, when applied to prime numbers $n$ of essentially the same size. In Section 5 we weaken the hypothesis given in the two previous sections and present a test for this larger family of numbers and some applications. In Section 6 we compare our algorithms with the AKS algorithm (when such comparison is valid), and we indicate some possible paths for future investigations. We include an explicit plausible conjecture.

This paper is modelled after [3]. The structure is very similar. The results on the theory of finite fields required can be found in many basic textbooks on finite fields or number theory, for example [11].

## 2 Preliminaries and Notation

Throughout the section $p$ denotes an odd prime number. Let $a$ be an integer coprime with $p$. The Legendre symbol $(\frac{a}{p})$ is defined by the formula

$$\left(\frac{a}{p}\right) = \begin{cases} 1 & : \text{ if there is an integer } x \text{ such that } x^2 \equiv 1 (\bmod p) \\ -1 & : \text{ otherwise .} \end{cases}$$

This symbol has the following properties:

1. If $ab$ is coprime with $p$ then $(\frac{ab}{p}) = (\frac{a}{p})(\frac{b}{p})$.

2. $(\frac{a}{p}) \equiv a^{\frac{p-1}{2}} \pmod{p}$.

The Legendre symbol can be extended multiplicatively to the Jacobi symbol replacing $p$ by an odd number $m$. That is, if $m = p_1...p_k$ and $(a, m) = 1$ then $(\frac{a}{m}) = (\frac{a}{p_1})....(\frac{a}{p_k})$. The Jacobi symbol also satisfies property (1) of the Legendre symbol above. Most important, it satisfies the well-known quadratic reciprocity law which we now state.

Let $m, n$ be odd and coprime numbers. Then,

1. $(\frac{-1}{n}) = (-1)^{\frac{n-1}{2}}$

2. $(\frac{2}{n}) = (-1)^{\frac{n^2-1}{2}}$.

3. $(\frac{m}{n}) = (\frac{n}{m})(-1)^{\frac{m-1}{2}\frac{n-1}{2}}$.

The proof of the quadratic reciprocity law can be found in most text books in number theory. As a reference we give [11].

Let $F_p$ denote the finite field with $p$ elements. For the sake of readability we recall some basic facts about the theory of finite fields that we shall employ below. These facts can also be found in many text books in the subject. We give [14] as a reference.

**Proposition 2.1** *Let $K$ and $E$ be finite fields containing $F_p$. Let $q = |K|$ and suppose $K \subseteq E$. Then,*

1. *$E$ has $q^d$ elements for some positive integer d.*

2. *$E$ is vector space of dimension d over $K$ $(d = [E : K])$.*

3. *$E$ is isomorphic to $K[x]/h(x)$ polynomial $h(x) \in K[x]$ of degree d irreducible over $K$.*

4. Let $C$ an algebraic closure of $F_p$ containing $E$. Then $E$ is the unique subfield of $C$ of dimension $d$ over $K$. It is the smallest subfield of $C$ containing a root $\theta$ of $h(x)$.

$$E = K(\theta) = \{f(\theta) | f(x) \in K[x], \ degree \ f(x) < d\}.$$

5. The multiplicative group $E^*$ of $K$ is cyclic of order $q^d - 1$.

6. The Galois group $G = Gal(E/K)$, that is, the group of automorphism of $E$ over $K$, is cyclic of order $d$. It is generated by the Frobenius automorphism $\sigma_q$ defined by $\sigma_q(\alpha) = \alpha^q$ for all $\alpha \in E$.

7. $[E : F_p] = [E : K][K : F_p]$.

Now let $K$ be a finite extension of $F_p$ with $q = |K|$. Let $K^\star$ be the multiplicative group and $g$ a generator of $K^\star$.

**Lemma 2.1** *For an element $\alpha$ of $K$, the following are equivalent*

1. $x^2 - \alpha^l$ is irreducible over $K$ for every odd integer $l$.

2. $x^2 - \alpha$ is irreducible over $K$.

3. $\alpha = g^t$ for some odd integer $t$.

4. $\alpha^{\frac{q-1}{2}} = -1$.

**Proof**

(1) $\Rightarrow$ (2) is trivial. Now let us prove (2) $\Rightarrow$ (3). Since $g$ is a generator, then $\alpha = g^t$ for some $t$. If $t = 2m$ then $x^2 - \alpha = x^2 - g^{2m} = (x - g^m)(x + g^m)$ is reducible. (3) $\Rightarrow$ (4) is obtained by noticing that $g^{\frac{q-1}{2}} = -1$ since $g$ is a generator. Hence, $\alpha^{\frac{q-1}{2}} = (-1)^t = -1$. Finally, to show (4) $\Rightarrow$ (1) suppose $x^2 - \alpha$ is reducible. Then, there is $\beta \in K$ such that $\beta^2 = \alpha$. So $\alpha^{\frac{q-1}{2}} = \beta^{q-1} = 1$ which contradicts the assumption. □

**Lemma 2.2** *Let $q = |K|$. Assume $q \equiv 1 \pmod 4$. If $x^2 - a$ is irreducible over $K$ and $\theta$ is a root of $x^2 - a$, then $x^2 - \theta$ is irreducible over $K(\theta)$.*

**Proof**

Note that $|K(\theta)| = q^2$. By Lemma 2.1 it is enough to prove that $\theta^{\frac{q^2-1}{2}} = -1$. Note that since $q \equiv 1 \pmod 4$ then $\frac{q+1}{2} = t$ is odd. Also, since $x^2 - a$ is irreducible over $K$, then $a^{\frac{q-1}{2}} = -1$. Hence,

$$\theta^{\frac{q^2-1}{2}} = ((\theta^2)^{\frac{q+1}{2}})^{\frac{q-1}{2}} = (-1)^t = -1.$$

□

**Corollary 2.1** *If $|K| = q \equiv 1 \pmod 4$ and $a \in K$ is such that $a^{\frac{q-1}{2}} = -1$, then the polynomial $x^{2^s} - a$ is irreducible over $K$ for all $s \geq 1$.*

**Proof**

Proceed inductively on $s$. Use Lemma 2.2 and part 7 of Proposition 2.1 □

We can now establish the following proposition

**Proposition 2.2**  1. If $p \equiv 1$ (mod 4) and $(\frac{a}{p}) = -1$, then $x^{2^s} - a$ is irreducible over $F_p$.

2. If $p \equiv 3$ (mod 4) and $(\frac{a}{p}) = (\frac{1-a}{p}) = -1$, then $x^{2^s} - 2x^{2^{s-1}} + a$ is irreducible over $F_p$.

**Proof**

The assertion (1) is a particular case of the Corollary 2.1 since $(\frac{a}{p}) \equiv a^{\frac{p-1}{2}} \equiv -1$ (mod $p$). In order to prove (2) let $\theta_1 = 1 + \sqrt{1-a}$. Since $(\frac{1-a}{p}) = -1$ then $F_p(\theta_1)$ has degree 2 over $F_p$. Hence it has $p^2 = q$ elements, so $q \equiv 1$ (mod 4). Moreover,

$$\theta_1^{\frac{p^2-1}{2}} = (\theta_1^{p+1})^{\frac{p-1}{2}} = ((1+\sqrt{1-a})(1-\sqrt{1-a}))^{\frac{p-1}{2}} = a^{\frac{p-1}{2}} = -1.$$

Corollary 2.1 implies that $x^{2^{s-1}} - \theta_1$ is irreducible over $F_p(\theta_1)$. A root $\theta$ of this polynomial satisfies,

$$(x^{2^{s-1}} - \theta_1)(x^{2^{s-1}} - \theta_1^p) = x^{2^s} - 2x^{2^{s-1}} + a$$

which belongs to $F_p[x]$. By part 7 of Proposition 2.1 it must be irreducible over $F_p$. □

## 3  Algorithm for the case $n \equiv 1 \pmod{4}$

Throughout this section we assume that $n \equiv 1$ (mod 4). Let $k = \nu_2(n-1)$. So $k \geq 2$. Let $a$ be an integer such that $(\frac{a}{n}) = -1$. Note for example that if $n = h \, 2^k + 1$ and $h \not\equiv 0$ (mod 3) then $n$ is either a multiple of 3 or $(\frac{3}{n}) = -1$. This is easily deduced from the quadratic reciprocity law. It follows that the algorithm that we will present in this section is deterministic for numbers of that form. Finally let $s = \lceil 2 \log \log n \rceil$. Note that $(\log n)^2 < 2^s < 2(\log n)^2$. We now describe the proposed Algorithm.

**Algorithm 1**

**Input** $n, a$: $n \equiv 1$ *(mod 4)*, $(\frac{a}{n}) = -1$.
  Let $k = \nu_2(n-1)$, $s = \lceil 2 \log \log n \rceil$.

1. *Verify properties of the Legendre Symbol and Proth's Theorem.*

    (a) *Let $A = a^{\frac{n-1}{2^k}}$. If $A^{2^{k-1}} \not\equiv -1$ (mod $n$), output* **composite**.

    (b) *If $k > (1/2) \log n$, output* **prime**.

2. *Verify $n$ is not a perfect power.*

    *If $n = d^e$ for some positive integers $d$ and $e$ with $e > 1$, output* **composite**.

3. Generate a set $\mathcal{S}$ of cardinality $2^{\max(s-k,0)}$

    $m = 1, \mathcal{S} = \{1\}, \mathcal{S}' = \{1\}$.

    While $(|\mathcal{S}| < 2^{\max(s-k,0)})\{$

        While $(m^{2^k} \pmod{n} \in \mathcal{S}')\{$

            $m \leftarrow m + 1$

        $\}$.

        If $m > |\mathcal{S}|2^k + 1$, output **composite**.

        If $(m,n) > 1$, output **composite**.

        If $(m^{2^k} - s', n) > 1$ for some $s' \in \mathcal{S}'$, output **composite**.

        $\mathcal{S} \leftarrow \mathcal{S} \bigcup \{m\}$.

        $\mathcal{S}' \leftarrow \mathcal{S}' \bigcup \{m^{2^k} \pmod{n}\}$

    $\}$.

4. For all $m \in \mathcal{S}$.

    If $(1+mx)^n \not\equiv (1+mx^n) \bmod (n, x^{2^s} - a)$, output **composite**.

    Output **prime**.

**Theorem 3.1** *The algorithm above returns prime if and only if $n$ is prime (as long as $n > 100$).*

**Theorem 3.2** *The running time of the algorithm is $\tilde{O}(2^{-\min(s,k)}(\log n)^6)$. Note that this is $\tilde{O}((\log n)^6)$ if $k = 2$ and is $\tilde{O}((\log n)^4)$ if $k \geq s$.*

The rest of this section is devoted to the proof of these theorems. We do this in a similar way as is done in [3], through the proofs of a series of lemmas.

**Lemma 3.1** *If $n$ is prime ($n > 100$), the algorithm returns prime.*

**Proof**

    Step (1a) of the algorithm can not return composite because of property 2 of the Legendre symbol.

    Step (2) can not return composite because $n$ is not a perfect power.

    Next we show Step (3) does not return composite. First note that if $k \geq s$ then the algorithm does not enter the first *while* loop, hence Step (3) cannot return composite in this case. So we may assume $k < s$. In this case the algorithm generates the set $\mathcal{S}$, that is, a sequence of integers $m_i$ with $i = 1, ..., 2^{s-k}$. $m_1 = 1$. Since $n$ is prime, the number of solutions of $x^{2^k} = 1$ in $F_n$ is at most $2^k$ (in fact it is exactly $2^k$ since the distinct powers of $A$ are solutions of this equation). It follows that $m_2 \leq 2^k + 1$. Inductively, using this same reasoning we deduce that $m_t \leq (t-1)2^k + 1$. Note that $t-1$ is the cardinality of the set $\mathcal{S}$ at that stage of the algorithm. It follows that under the assumption that $n$ is prime, $m > |\mathcal{S}|2^k + 1$ cannot occur.

It also follows that each $m_i \leq (2^{s-k} - 1)2^k + 1 < 2^s < 2(\log n)^2 < n$ (this last inequality certainly occurs if $n > 100$). Hence, in the algorithm $(m, n) > 1$ cannot occur. Finally, since $m_i^{2^k} \not\equiv m_j^{2^k} \pmod{n}$ for all $j < i$, then $(m^{2^k} - s', n) > 1$ cannot occur. This concludes the analysis for Step (3).

Since $(1 + mx)^n \equiv 1 + mx^n \pmod{n}$ then $(1 + mx)^n \equiv 1 + mx^n \pmod{(n, x^{2^s} - a)}$, so Step (4) does not return composite. □

We assume from now on that the output of the algorithm is prime.

**Lemma 3.2** *Suppose that the algorithm has passed step (1a), that is, it has verified $A^{2^{k-1}} \equiv -1 \pmod{n}$. Then, we have*

1. *$\nu_2(d-1) \geq k$ for all divisors of $n$.*

2. *There is a prime divisor $p$ of $m$ for which $\nu_2(p-1) = k$. For such prime $p$, $\left(\frac{a}{p}\right) = -1$.*

**Proof**

1. It is enough to prove it for prime divisors $d$ of $n$. The hypothesis implies $A^{2^{k-1}} \equiv -1 \pmod{d}$, whence $\text{ord}_d(A) = 2^k$, so $\nu_2(d-1) \geq k$.

2. If every prime divisor $q$ of $n$ were to satisfy $\nu_2(q-1) > k$, then so would the product, that is, $n$. Let $p$ a prime divisor of $n$ satisfying $\nu_2(p-1) = k = \nu_2(n-1)$. Let $t = \frac{p-1}{2^k}$. Note that $t$ is odd. Hence

$$\left(\frac{A}{p}\right) \equiv A^{\frac{p-1}{2}} \equiv (A^t)^{2^{k-1}} \equiv (-1)^t \equiv -1 \pmod{p}.$$

Since $A = a^{\frac{n-1}{2^k}}$ and $\frac{n-1}{2^k}$ is odd, then we get the result. □

**Lemma 3.3** *If the algorithm output prime at Step (1b) then $n$ is prime.*

**Proof**
This follows Proth's Theorem [13]. Let us recall its statement: if $\nu_2(n-1) > (1/2)\log n$ and $\left(\frac{a}{n}\right) = -1$, then $n$ is prime if and only if $a^{\frac{n-1}{2}} \equiv -1 \pmod{n}$. □

Now we assume $n$ has passed Step (1b) (so $k \leq 1/2 \log n$). We let $p$ be a prime divisor of $n$ satisfying $\nu_2(p-1) = k = \nu_2(n-1)$. Since $\left(\frac{a}{p}\right) = -1$, then by Proposition 2.2, the polynomial $x^{2^s} - a$ is irreducible over $F_p$. Let $\theta$ be a root of the polynomial in an algebraic closure $\mathcal{C}$ of $F_p$, let $K = F_p(\theta)$ and $K^\star$ its multiplicative group. Every $\alpha \in K^\star$ is $\alpha = f(\theta)$ for some (unique) non-zero polynomial $f(x) \in F_p[x]$ of degree $t < 2^s$. Let $m$ be an integer. We denote by $r_m$ the multiplicative homomorphism of $K^\star$ consisting in raising to the $m$-th power. We denote by $\sigma_m$ the linear map of $K$ defined by $\sigma_m(\alpha) = f(\theta^m)$, where $f(x)$ is the unique polynomial mentioned above.

**Lemma 3.4** *For an integer $m$ the following are equivalent:*

1. $\theta^m$ is a root of $\mathrm{irr}_\theta(x) = x^{2^s} - a$.

2. $a^m = a$ (in $F_p$).

3. $\sigma_m(h(\theta)) = h(\theta^m)$ for all $h(x) \in F_p[x]$.

4. $\sigma_m \in Gal(K/F_p)$.

**Proof**

That (1) implies (2) is clear since $(\theta^m)^{2^s} = a^m$. To see that (2) implies (3) write $h(x) = f(x) + (x^{2^s} - a)p(x)$ where $f(x)$ has degree less than $2^s$. By definition of $\sigma_m$ we have

$$\sigma_m(h(\theta)) = f(\theta^m) = h(\theta^m) - (a^m - a)p(\theta^m) = h(\theta^m).$$

To prove that (3) implies (4) note that since $\sigma_m$ is clearly a linear map over $F_p$ we only have to prove it is multiplicative, and this is trivial. Finally, (4) implies (1) is also evident: just note that $\sigma_m(\theta) = \theta^m$ is a conjugate of $\theta$ over $F_p$, hence, it must be a root of $irr_\theta(x)$.

□

In particular, since $a^n \equiv a \pmod{n}$, this lemma implies that $\sigma_n \in Gal(K/F_p)$, so it must be a power of the Frobenius automorphism $\sigma_p^i = \sigma_{p^i}$. The idea will be to show that, under certain conditions that are met if the algorithm outputs prime in its last step, this implies that $n = p^i$. We still need quite a few observations before reaching that conclusion.

Write $n = p^l d$. Then, from $a^n = a$ and $a^{p^l} = a$ it is easy to deduce that $a^d = a$. So $\sigma_d$ is also an automorphism. Moreover, so is $\sigma_{d^i p^j}$ for all $i, j \geq 0$. More generally if $m_1$ and $m_2$ satisfy the equivalent conditions of the previous lemma then so does $m_1 m_2$ and it is also easy to verify that $\sigma_{m_1 m_2} = \sigma_{m_1} \circ \sigma_{m_2}$. Similarly, if $m_1$ and $m_1 m_2$ satisfy these conditions, then so does $m_2$. On the other hand, if $m$ satisfies any of the equivalent conditions of the previous lemma then the product $\sigma_m r_{-m}$ is also a multiplicative homomorphism of $K^\star$ since it is a product of homomorphisms. It follows that

$$G_m = \mathrm{Ker}\sigma_m r_{-m} = \{f(\theta) \in K^\star | f(\theta^m) = f(\theta)^m\}$$

is a subgroup of $K^\star$, hence cyclic generated by, say, $g_m(\theta)$. We now study the properties of these cyclic groups.

**Lemma 3.5** *Suppose $m_1$ and $m_2$ satisfy any of the equivalent conditions of lemma 3.4, then:*

1. *For all $i \geq 0$, $G_{p^i} = K^\star$.*

2. *$G_{m_1} \cap G_{m_2} \subseteq G_{m_1 m_2}$.*

3. *$|G_{m_i}|$ divides $m_i^{2^s} - 1$. In particular $(m_i, |G_{m_i}|) = 1$.*

4. $G_{m_1 m_2} \cap G_{m_1} \subseteq G_{m_2}$.

**Proof**

1. That $G_1 = K^\star$ is trivial. Let $\alpha \in K^\star$ then $\sigma_{p^i}(\alpha) = \sigma_p^i(\alpha) = \alpha^{p^i}$, since $\sigma_p$ is the Frobenious automorphism.

2. Let $\alpha \in G_{m_1} \cap G_{m_2}$. Then, $\sigma_{m_1}(\alpha) = \alpha^{m_1}$ and $\sigma_{m_2}(\alpha) = \alpha^{m_2}$. It follows that

$$\sigma_{m_1 m_2}(\alpha) = \sigma_{m_1}(\sigma_{m_2}(\alpha)) = \sigma_{m_1}(\alpha^{m_2}) = (\sigma_{m_1}(\alpha))^{m_2} = (\alpha^{m_1})^{m_2} = \alpha^{m_1 m_2}.$$

This implies $\alpha \in G_{m_1 m_2}$.

3. Let $\alpha$ be a generator of $G_{m_i}$. By part 2 of this lemma $\alpha$ belongs to $G_{m_i^{2^s}}$. On the other hand, since $\sigma_{m_i}$ is an automorphism of $K$ then $\sigma_{m_i}^{2^s} = $ identity. So $\alpha^{m_i^{2^s}} = \sigma_{m_i^{2^s}}(\alpha) = \sigma_{m_i}^{2^s}(\alpha) = id(\alpha) = \alpha$. So $\alpha^{m_i^{2^s}-1} = 1$. Hence $|G_{m_i}| = \mathrm{ord}(\alpha)$ divides $m_i^{2^s} - 1$. In particular $(m_i, |G_{m_i}|) = 1$.

4. Let $\alpha \in G_{m_1 m_2} \cap G_{m_1}$. Then

$$(\alpha^{m_2})^{m_1} = \alpha^{m_1 m_2} = \sigma_{m_2}(\sigma_{m_1}(\alpha)) = \sigma_{m_2}(\alpha^{m_1}) = (\sigma_{m_2}(\alpha))^{m_1}$$

so $(\alpha^{m_2})^{m_1} = (\sigma_{m_2}(\alpha))^{m_1}$. By the previous item of this lemma there is an integer $t$ such that $tm_1 \equiv 1 \pmod{|G_{m_1}|}$. Raising to this $t$ we obtain $(\alpha^{m_2})^{m_1 t} = (\sigma_{m_2}(\alpha))^{m_1 t}$. Note that $\sigma_{m_2}(\alpha)$ has the same order than $\alpha$. Hence $\alpha^{m_2} = \sigma_{m_2}(\alpha)$. □

Write $n = p^l d$ where $d$ is coprime with $p$. We use the previous lemma to obtain the following result.

**Corollary 3.1** *For all $i, j \geq 1$, $G_n \subseteq G_{p^i} G_{d^j}$.*

**Proof**

$$G_n = G_{dp^l} = G_{dp^l} \cap G_{p^l} \subseteq G_d \subseteq G_{d^i} = G_{d^i} \cap G_{p^j} \subseteq G_{d^i p^j}.$$ □

**Corollary 3.2** *(Analogous to Lemma 4.6 in [3]). If $m_1$ and $m_2$ satisfy any of the equivalent conditions of lemma 3.4, then $\sigma_{m_1} = \sigma_{m_2}$ implies $|G_{m_1} \cap G_{m_2}|$ divides $m_1 - m_2$.*

**Proof**

Let $\alpha \in G_{m_1} \cap G_{m_2}$. Then $\alpha^{m_1} = \sigma_{m_1}(\alpha) = \sigma_{m_2}(\alpha) = \alpha^{m_2}$, thus $\alpha^{m_1 - m_2} = 1$. Since, $G_{m_1} \cap G_{m_2}$ is a cyclic group, then $|G_{m_1} \cap G_{m_2}|$ divides $m_1 - m_2$. □

The following lemma is very important because it shows how to obtain $2^{\min(k,s)}$ monomials in $G_n$ from one iteration in Step 4 of the algorithm. This is the reason why the complexity of the algorithm improves as $k$ grows.

**Lemma 3.6**  1. *Suppose $k < s$. If for some integer $m$, we have $(1 + m\theta) \in G_n$, then $(1 + mA^i\theta) \in G_m$ for $i = 1, 2, ..., 2^k$.*

2. *Suppose $k \geq s$. Let $B = A^{2^{k-s}}$. If $(1+\theta) \in G_n$, then $(1+B^i\theta) \in G_n$ for $i = 1, 2, ..., 2^s$.*

**Proof**

1. Recall $G_n$ is a group, so $(1+m\theta) \in G_n$ implies $(1+m\theta)^{p^i} = (1+m\theta^{p^i}) \in G_n$. The elements $\theta^{p^i}$ are the Galois conjugates of $\theta$ in $F_p[\theta]$. Since $\theta^{2^s} = A$, then the conjugates are of the form $\theta\zeta$, where $\zeta^{2^s} = 1$. Since $k \leq s$ every $A^i$ satisfies $(A^i)^{2^s} = 1$. So the $A^i$ are among the possible values for $\zeta$. In particular, $(1 + mA^i\theta) \in G_n$.

2. Same as in (1) by noting that $B$ is a primitive $2^s$-th root of 1 in $F_p$  $\square$

**Lemma 3.7** *If the algorithm ouputs prime at Step 4, then $|G_n| \geq 2^{2^s}$.*

**Proof**

Assume first that $k < s$. Again we denote by $m_i$, with $i = 1, \ldots, 2^{s-k}$, the sequence of elements of the set $S$ generated by the algorithm in Step (3). We claim that $m_i A^j$ for $i = 1, 2, \ldots, 2^{s-k}$ and $j = 1, \ldots, 2^k$ are all different and non-zero in $F_p$. To see this recall that $A$ has order $2^k$ in $F_p$. Hence $A^j$ is non zero in $F_p$ for all $j$ and they are all different for $j = 1, \ldots, 2^k$. The algorithm verifies $(m_i, n) = 1$. Hence, the $m_i A^j$ are all non zero in $F_p$. Assume $m_i A^j = m_{i'} A^{j'}$ in $F_p$. Raising to the $2^k$th power we get $m_i^{2^k} = m_{i'}^{2^k}$ in $F_p$, but since the algorithm verified that $m_i^{2^k} - m_{i'}^{2^k}$ is coprime with $n$, then we must have $i = i'$ whence we deduce that $j = j'$. So we have $2^s$ different non-zero elements of $F_p$. Denote them by $t_r$ for $r = 1, \ldots, 2^s$. The algorithm verifies that $(1 + m_i\theta) \in G_n$ for each $i = 1, 2, \ldots, 2^{s-k}$. It follows from the previous lemma that $(1 + t_r\theta) \in G_n$ for $r = 1, 2, \ldots, 2^s$.

If, on the other hand, $k > s$, then the algorithm verifies that $(1+\theta) \in G_n$, and, again, using the previous lemma, we get $(1 + B^r\theta) \in G_n$ for $r = 1, 2, \ldots, 2^s$. So in both cases we obtain $2^s$ different monomials in $G_n$. To simplify we always denote these $(1 + t_r\theta) \in G_n$ for $r = 1, 2, \ldots, 2^s$. Since $G_n$ is a group it contains the set $\mathcal{T}$ defined as

$$\mathcal{T} = \left\{ \prod_{r=1}^{2^s}(1 + t_r\theta)^{\epsilon_r} \mid \epsilon_r \in \mathbb{Z}^+, \sum \epsilon_r < 2^s \right\}.$$

Every element of $\mathcal{T}$ is of the form $f(\theta)$ for some $f(x)$ of degree less than $2^s$. Since all $t_r$ are different in $F_p$ then the polynomials $f(x)$ corresponding to different choices of $\epsilon_i$ are different in $F_p[x]$. Since the degrees are less than $2^s$, then the corresponding elements of $\mathcal{S}$ are different.

$\mathcal{T}$ contains properly the set

$$\mathcal{T}_1 = \left\{ \prod_{r=1}^{2^s}(1 + t_r\theta)^{\epsilon_r} \mid \epsilon_r \in \{0,1\}, \sum \epsilon_r < 2^s \right\}$$

with cardinality $2^{2^s} - 1$. Hence, $\mathcal{T}$ has at least $2^{2^s}$. Therefore $|G_n| \geq 2^{2^s}$. □

We are now ready to complete the proof of Theorem 3.1.

**Proof of Theorem 3.1**

It remains to prove that if the algorithm outputs prime in the last step, then $n$ is prime. Assume $n$ has more than one prime divisor. Hence, $n = p^l d$ where $(d, p) = 1$ and $d > 1$. We know that $\sigma_{p^i d^j} \in Gal(K/F_p)$ for all $i, j \geq 0$. Since $Gal(K/F_p)$ has order $2^s$ it follows from the pigeon hole principle that there exist two different pairs $(i_1, j_1)$ and $(i_2, j_2)$ with $0 \leq i_1, j_1, i_2, j_2 \leq [\sqrt{2^s}]$ such that $\sigma_{p^{i_1} d^{j_1}} = \sigma_{p^{i_2} d^{j_2}}$. It follows from Corollary 3.2 that

$$|G_{p^{i_1} d^{j_1}} \cap G_{p^{i_2} d^{j_2}}| \text{ divides } p^{i_1} d^{j_1} - p^{i_2} d^{j_2}.$$

Hence, from Corollary 3.1 we obtain

$$|G_n| \text{ divides } p^{i_1} d^{j_1} - p^{i_2} d^{j_2}.$$

Note that $p^{i_1} d^{j_1} - p^{i_2} d^{j_2} < n^{[\sqrt{2^s}]} \leq n^{\sqrt{2^s}}$. Also note that from $s = \lceil 2 \log \log n \rceil$ one can easily deduce that $2^{2^s} > n^{\sqrt{2^s}}$. It follows from Lemma 3.7 that $|G_n| > n^{\sqrt{2^s}}$. So we obtain $p^{i_1} d^{j_1} = p^{i_2} d^{j_2}$. But this is not possible because $p$ and $d$ are coprime and $(i_1, j_1) \neq (i_2, j_2)$. Hence $d = 1$. So, $n = p^l$. Since $n$ passed Step 2 of the algorithm ($n$ is not a perfect power) we conclude $l = 1$, so $n = p$ □

**Analysis of Complexity. Proof of Theorem 3.2**

Step 1 involves the calculation of $a^{\frac{n-1}{2}} \pmod{n}$ which takes $\tilde{O}((\log n)^2)$ time using the fast Fourier transform.

Step 2, as in [3] takes $\tilde{O}((\log n)^3)$.

Step 3. If $k \geq s$ the algorithm does not enter the *while* loop, so in this case this step has no cost. When $k < s$, every integer $m$ that the algorithm deals with is less than $2^s$. For each of these integers $m$, it computes $m^{2^k} \pmod{n}$. It follows that the algorithm calculates $m^{2^k}$ for at most $2^s$ different values of $m$ (in practice much less than this). This involves $k 2^s \leq s 2^s$ modular multiplications (multiplications mod $n$). Using the fast Fourier transform these computations take at most $\tilde{O}((\log n)^3)$. On the other hand, the algorithm in this Step computes less than $2^{2(s-k)}$ gcd's. This takes $2^{2(s-k)} \tilde{O}((\log n)) = 2^{-2k} \tilde{O}((\log n)^5)$ time.

Step 4: This is the part of the computation that will determine the complexity of the algorithm. It involves $2^{\max(s-k,0)}$ iterations, where by iteration we mean the computation of $(1 + m_i x)^n \mod (n, x^{2^s} - a)$. Using fast exponentiation each iteration takes at most $2 \log n$ multiplications in the field $K$. Using the fast

Fourier transform each of these involves $O(2^s s)$ modular multiplications, and likewise each of these takes $\tilde{O}(\log n)$ time. We must add that the reduction modulo $x^{2^s} - a$ is necessary after multiplications of elements in $K$, but these are done with $2^s$ modular multiplications, which does not affect complexity. So each iteration takes $\tilde{O}((\log n)^4)$. Hence this step takes

$$2^{\max(s-k,0)} \, \tilde{O}((\log n)^4) = 2^{-\min(s,k)} \, \tilde{O}((\log n)^6),$$

and so does the algorithm. □

## 4 Algorithm for $n \equiv -1 \pmod{4}$

Throughout this section we assume that $n \equiv -1 \pmod 4$, and $k = \nu_2(n+1)$. In particular $k \geq 2$. We assume that an integer $a$ is given such that $(\frac{a}{n}) = (\frac{1-a}{n}) = -1$. Note for example that if $n = h2^k - 1$ and $h \not\equiv 0 \pmod 3$ then $n$ is either a multiple of 3 or $(\frac{3}{n}) = (\frac{1-3}{n}) = -1$. This is easily deduced from the quadratic reciprocity law. It follows that the algorithm presented in this section is deterministic for numbers of that form. Further we let $t = \lceil 2 \log \log n \rceil + 1$, noting that $t = s + 1$. Hence $2 \, (\log n)^2 < 2^t < 4 \, (\log n)^2$. We now describe the proposed Algorithm.

**Algorithm 2**
**Input** $n, a$: $n \equiv -1$ (mod 4), $(\frac{a}{n}) = (\frac{1-a}{n}) = -1$.

Compute $k = \nu_2(n+1)$, $t = \lceil 2 \log \log n \rceil + 1$.

1. Verify properties of the Legendre Symbol, the Frobenius automorphism and Lucas-type Theorem.

   (a) If $a^{\frac{n-1}{2}} \not\equiv -1 \pmod{n}$, output **composite**.

   (b) If $(1 + \sqrt{1-a})^n \not\equiv 1 - \sqrt{1-a} \pmod{n}$ output **composite**.

   (c) If $k > 1/2 \log n$ output **prime**.

2. Verify $n$ is not a perfect power.

   If $n = d^e$ for some positive integer $e$, output **composite**.

3. Finding a sequence of $m_i$'s.

   For $m = 1$ to $2^{\max(t-k,0)}$

   $$\text{If } (m,n) > 1, \text{ output } \textbf{composite}.$$

4. Finding elements in $G_n$.

   For $m = 1$ to $2^{\max(t-k-1,0)}$

$$\text{If } (1 + mx)^n \not\equiv (1 + mx^n) \pmod{n, x^{2^{t+1}} - 2x^{2^t} + a} \text{ output } \textbf{composite}.$$

output **prime**

**Theorem 4.1** *The algorithm above returns prime if and only if $n$ is prime (assuming $n > 25$).*

**Theorem 4.2** *The running time of the algorithm is $\tilde{O}(2^{-\min(s,k)}(\log n)^6)$.*

The proofs of these results are analogous to the theorems in the previous section. However, in many occasions, the analogy is not immediate. In these cases, we will go over the necessary lemmas and give detailed proofs.

**Lemma 4.1** *If $n$ is prime, the algorithm returns prime.*

**Proof**

Step 1 cannot output composite: in the first place because of the properties of the Legendre Symbol, and secondly because of the properties of the Frobenius automorphism. The rest proceeds as in the case $n \equiv 1 \pmod{4}$, except that in Step 3 we only need $n > 25$ to make sure that $2^{\max(t-k,0)} < n$. □

We assume now that the output of the algorithm is prime.

**Lemma 4.2** *Let $n, a, 1 - a$ as in the input of the algorithm, and $k = \nu_2(n+1)$. Suppose $a^{\frac{n-1}{2}} \equiv -1 \pmod{n}$ and that $(1 + \sqrt{1-a})^n \equiv 1 - \sqrt{1-a} \pmod{n}$. Then,*

1. *Every prime divisor $q$ of $n$ satisfies either*

   (a) $q \equiv 1 \pmod{2^{k+1}}$ *or*

   (b) $q \equiv -1 \pmod{2^k}$

   *It satisfies (a) if and only if $(\frac{1-a}{q}) = (\frac{a}{q}) = 1$.*
   *It satisfies (b) if and only if $(\frac{1-a}{q}) = (\frac{a}{q}) = -1$.*

2. *There exists a prime divisor $p$ of $n$ such that $\nu_2(p+1) = \nu_2(n+1) = k$. For such $p$, $(\frac{1-a}{p}) = (\frac{a}{p}) = -1$.*

**Proof**

1. Let $q$ be a prime divisor of $n$. We first note that $(\frac{a}{q}) = 1$ if and only if $q \equiv 1 \pmod{4}$. Recall $\frac{n-1}{2}$ is odd. Hence $(\frac{a}{q}) = (\frac{a}{q})^{\frac{n-1}{2}} = (\frac{-1}{q}) = (-1)^{\frac{q-1}{2}}$.

   Next we show that $(1 + \sqrt{1-a})^{\frac{n^2-1}{2}} \equiv -1 \pmod{n}$. This is true since

   $$(1 + \sqrt{1-a})^{\frac{n^2-1}{2}} = ((1 + \sqrt{1-a})^{n+1})^{\frac{n-1}{2}} = ((1 - \sqrt{1-a})(1 + \sqrt{1-a}))^{\frac{n-1}{2}} = a^{\frac{n-1}{2}} \equiv -1 \pmod{n}.$$

   Now suppose $(\frac{1-a}{q}) = 1$. Then, $F_q(\sqrt{1-a}) = F_q$. Since $(1 + \sqrt{1-a})^{\frac{n^2-1}{2}} \equiv -1 \pmod{n}$ then $(1 + \sqrt{1-a})^{\frac{n^2-1}{2}} = -1$ in $F_q$. But $\nu_2(n+1) = k$ implies $\nu_2(n^2 - 1) = k + 1$. So the element $(1+\sqrt{1-a})^{\frac{n^2-1}{2^{k+1}}}$ has order $2^{k+1}$ in $F_q$ so $q \equiv 1 \pmod{2^{k+1}}$. In particular, $(\frac{a}{q}) = 1$ by our first remark.

Suppose now that $(\frac{1-a}{q}) = -1$. Then $F_q(\sqrt{1-a}) = F$ has $q^2$ elements. Again, $(1+\sqrt{1-a})^{\frac{n^2-1}{2}} = -1$ in $F$, so

$$(1+\sqrt{1-a})^{\frac{n^2-1}{2}} = \left((1+\sqrt{1-a})^{n-1}\right)^{\frac{n+1}{2}} = \left(\frac{1-\sqrt{1-a}}{1+\sqrt{1+a}}\right)^{\frac{n+1}{2}} = -1.$$

Note that in $F_q(\sqrt{1-a})$, the element $\beta = (\frac{1-\sqrt{1-a}}{1+\sqrt{1+a}}) = (1+\sqrt{1-a})^{q-1}$ lies in the unique subgroup of $F^\star$ of order $q+1$. $\beta^{\frac{n+1}{2^k}}$ has order $2^k$ so $q \equiv -1 \pmod{2^k}$. Also, $(\frac{a}{q}) = -1$ by our first remark.

2. Since $(\frac{1-a}{n}) = -1$ then there must be a prime divisor of $n$ such that $(\frac{1-a}{q}) = -1$. So $q \equiv -1 \pmod{2^k}$. If all primes satisfying $(\frac{1-a}{q}) = -1$ satisfy $q \equiv -1 \pmod{2^{k+1}}$, then by part 1, $n$ would satisfy $n \equiv \pm 1 \pmod{2^{k+1}}$. But, $\nu_2(n+1) = k$ implies this is not possible. So there is $p/n$ such that $\nu_2(p+1) = k$. For such $p$, which is congruent to $-1 \pmod 4$, we must have $(\frac{a}{p}) = -1$. Hence, we also must have $(\frac{1-a}{p}) = -1$ since we just proved that $(\frac{1-a}{p}) = 1$ implies $(\frac{a}{p}) = 1$ □

**Corollary 4.1** *If the algorithm outputs prime in Step 1c, then $n$ is prime.*

**Proof**

This is a small variation of the statement of a Lucas-type theorem. In any case, it is deduced easily from the previous lemma by noting that $k > 1/2 \log n$ is the same as $2^k > \sqrt{n}$, so the possible prime divisors are too large. □

Assume now that $n$ passed Step 1 of the algorithm and let $p$ the prime divisor of $n$ for which $\nu_2(p+1) = k$. We let $F = F_p(\sqrt{1-a})$ and $K = F_p(\theta)$ where $\theta$ is a root of the polynomial $x^{2^{t+1}} - 2x^{2^t} + a = \mathrm{irr}_\theta(x)$ which is irreducible by Proposition 2.2. We also note that $K = F(\theta)$ and $\theta$ is a root of $x^{2^t} - (1+\sqrt{1-a})$ or $x^{2^t} - (1-\sqrt{1-a})$, which are both irreducible over $F$. For simplicity we will assume $\theta$ is a root of the first of these two polynomials. The roots of the other one are also roots of $\mathrm{irr}_\theta(x)$. Let $\sigma_m$ defined as in the previous section by $\sigma_m(f(\theta)) = f(\theta^m)$ when $\deg f(x) < 2^s$. We need this lemma:

**Lemma 4.3** *For an integer $m$ the following are equivalent:*

1. $\theta^m$ is a root of $\mathrm{irr}_\theta(x) = x^{2^{t+1}} - 2x^{2^t} + a$.

2. $(1+\sqrt{1-a})^m = 1 \pm \sqrt{1-a}$ in $F$.

3. $\sigma_m(h(\theta)) = h(\theta^m)$ for all $h(x) \in F_p[x]$.

4. $\sigma_m \in \mathrm{Gal}(K/F_p)$.

We skip the proof as it is quite similar to that of lemma 3.4 of previous section.

When $\sigma_m$ is an automorphism we let

$$G_m = \{\alpha \in K : \sigma_m(\alpha) = \alpha^m\}.$$

Then $G_m$ is a cyclic subgroup of $K^\star$. Now write $n = p^l d$, where $p$ and $d$ are coprime. As in the previous section, we can use the above lemma to show that $\sigma_{d^i p^j} \in Gal(K/F_p)$ for all $i, j \geq 0$. Moreover we carry over Lemma 3.5, Corollary 3.1 and Corollary 3.2 in the new environment.

Let

$$\alpha = (1 + \sqrt{1-a})^{\frac{n^2-1}{2^{k+1}}}.$$

We have the following lemma analogous to Lemma 3.6.

**Lemma 4.4** Let $\beta = \alpha^{2^{\max(k+1-t,0)}}$. If $(1 + m\theta) \in G_n$ for some $m \neq 0$ in $F_p$, then $(1 + m\beta^i \theta) \in G_n$ for $i = 1, \ldots, 2^{\min(k+1,t)}$.

**Proof**

Proceed as in Lemma 3.7, since the conjugates of $\theta$ over the field $F$ are of the form $\theta \zeta$ where $\zeta^{2^t} = 1$. The powers of $\beta$ are among the latter. □

Next we estimate the size of $G_n$.

**Lemma 4.5** If the algorithm outputs prime in the last step then $|G_n| > 2^{2^t}$.

**Proof**

The algorithm verifies that every integer less than $2^{\max(t-k,0)}$ is coprime with $n$, hence they are all different and non-zero in $F_p$. Let $\gamma_{ij} = m_i \beta^j$ for $i = 1, 2, \ldots, 2^{\max(t-k-1,0)}$ and $j = 1, 2, \ldots, 2^{\min(k+1,t)}$. There are $2^t$ $\gamma_{ij}$'s. We claim they are all different and non-zero in $F$. Suppose $m_i \beta^j = m_{i'} \beta^{j'}$. Then $\frac{m_i}{m_{i'}} = \beta^{j'-j}$. Since the only powers of $\beta \in F_p$ are $2^{\min(k+1,t)}$ and $2^{\min(k,t-1)}$ (the other powers of $\beta$ are in $F - F_p$ we get: either $\beta^j = \beta^{j'}$, in which case $m_i = m_{i'}$ leading to $i = i'$, or $\beta^{j-j'} = -1$, in which case $m_i = -m_{i'}$. But then we have $m_i + m_{i'} = 0$ in $F_p$. Since $m_i + m_{i'} < 2^{\max(t-k,0)}$ and the algorithm verified in Step 3 that these were coprime with $n$ we get our claim. Next, since the algorithm verified that $(1 + m_i \theta) \in G_n$ for each $i$, it follows from the previous lemma that each of the $(1 + \gamma_{ij} \theta) \in G_n$. Therefore $G_n$ contains $2^t$ different monomials over $F$, and, as in the previous Section, we get the result. □

**Proof of Theorem 4.1**

Again this proceeds along the lines of the proof of Theorem 3.1. The only difference is that now $Gal(K/F_p)$ has order $2^{t+1}$ and $G_n$ has at least $2^{2^t}$ elements. The fact that $2^{2^t} > n^{\sqrt{2^{t+1}}}$ is easily derived from $2^{2^s} > n^{\sqrt{2^s}}$, keeping in mind that $t = s + 1$. □

**Remark 4.1** *Note that if $b, c$ are given integers such that $(\frac{b^2+c^2}{n}) = -1$ then $a = (bc^{-1})^2 + 1$ satisfies $(\frac{a}{n}) = (\frac{1-a}{n}) = -1$. This is easy to verify noting that $(\frac{-1}{n}) = -1$ since $n \equiv -1 \pmod{4}$. Alternatively, one could replace the polynomial in the algorithm by the polynomial $x^{2^{t+1}} - 2bx^{2^t} + (b^2 + c^2)$, which is also irreducible in $F_p$ under the assumption $(b^2 + c^2)^{\frac{n-1}{2}} \equiv -1 \pmod{n}$ and $(x + iy)^n \equiv (x - iy) \pmod{n}$.*

**Remark 4.2** *We note that the same polynomial used in the algorithm of this Section could have been used in the algorithm of the previous section, that is, for numbers $n \equiv 1 \pmod{4}$, with no additional hypothesis on $a$. To see this, notice that if $(\frac{a}{n}) = -1$ and $(\frac{1-a}{n}) = 1$ then $(\frac{a^{-1}}{n}) = -1$ and*

$$\left(\frac{1 - a^{-1}}{n}\right) = \left(\frac{-a^{-1}(1-a)}{n}\right) = \left(\frac{a^{-1}}{n}\right) = -1.$$

*So the pair $a, 1-a$ is achieved at most at the cost of computing $a^{-1}$. Hence, by Proposition 2.2 the polynomial $x^{2^{t+1}} - 2x^{2^t} - a$ is irreducible. However the algorithm we presented for numbers $n \equiv 1 \pmod{4}$ runs about four times faster than the other one. This is so, even though the number of iterations performed by both algorithms is the same, since the degree of the polynomial used in this Section is four times the degree of the polynomial used in the previous one.*

**Analysis of Complexity: Proof of Theorem 4.2**

The proof is similar to the proof of Theorem 3.2. We note that the cost of Step (3) is $\tilde{O}((\log n)^3)$, which is less than the cost of Step (3) of the Algorithm in the previous section, because the number of gcd's computed is much less in this algorithm. However, that does not improve complexity of the algorithm since it is Step (4) the Step that determines its complexity. In the previous remark we compared the speed of the two algorithms given in this paper. This comparison cannot be deduced from the notation used in the statement of the theorem, which is standard notation. □

## 5 Weakly Conditioned and Unconditioned tests

### 5.1 The case $n \equiv 1 \pmod{4}$

Let $n \equiv 1 \pmod{4}$. Let $k = \mu_2(n-1)$, so $k \geq 2$. This time we assume integers $a$ and $u$ are given, $1 \leq u \leq k$ such that $a^{\frac{n-1}{2^u}} \equiv -1 \pmod{n}$. Note that $u = 1$ is the case we dealt with in Section 3. At the other end, when $u = k$ there is always such an $a$, namely $a = -1$. Hence we will refer to the case $u = k$ as the *unconditional case*. We will produce a deterministic primality test for all such numbers. The complexity of the primality test we will give will depend also on $u$. The optimal performance occurs when $u = 1$ and the worst case is $u = k$.

We note that if $n = h2^k + 1$ is prime, and $h \not\equiv 0 \pmod{5}$ then either $5^{\frac{n-1}{2}} \equiv -1 \pmod{n}$ or $5^{\frac{n-1}{4}} \equiv -1 \pmod{n}$ or $n$ is a multiple of 5. This can be deduced from the law of biquadratic reciprocity. This fact was used in [6] to produce a deterministic primality test for numbers of that form provided $k > 1/2 \log n$.

Combining this observation with the one made at the beginning of Section 3. We deduce that every number of the form $n = h2^k + 1, h \not\equiv 0 \pmod{15}$ is either a multiple of 3 or 5 or can be tested using $a = 3$ or $a = 5$ and $u = 1$ or $u = 2$.

Again we let $s = \lceil 2 \log \log n \rceil$. We now present the algorithm in the form of a theorem.

**Theorem 5.1** Let $n \equiv 1 \pmod{4}$. Let $k = \nu_2(n-1)$. Let $s = \lceil 2 \log \log n \rceil$. Let $a$ and $u$ integers, $1 \leq u \leq k$ and such that $a^{\frac{n-1}{2^u}} \equiv -1 \pmod{n}$. Let $\mathcal{S}$ be a set of integers, $|\mathcal{S}| = 2^{\max(s-k+2(u-1),0)}$ such that for any pair $m, m'$ of different elements of $\mathcal{S}$, $(m^{2^{k+1-u}} - m'^{2^{k+1-u}}, n) = 1$ and such that every element of $\mathcal{S}$ is coprime with $m$. Suppose also that for every $m \in \mathcal{S}$ we have $(1 + mx)^n \equiv (1 + mx^n) \bmod (n, x^{2^{s+2(u-1)}} - a)$ and that $n$ is not a perfect power. Then, $n$ is prime.

**Proof** (Sketch)

Let $r = s + u - 1$. Let $f(x) = x^{2^{s+2(u-1)}} - a = x^{2^{r+u-1}} - a$. We enumerate some facts without a proof that can be deduced as in Section 3.

1. The equation $a^{\frac{n-1}{2^u}} \equiv -1 \pmod{n}$ implies that every prime divisor $q$ of $n$ satisfies $\nu_2(q-1) \geq k-u+1$.

2. There is a prime $p$ dividing $n$ such that $\nu_2(p - 1) \leq k$.

   Let $p$ be such a prime and $\theta$ a root of $f(x)$ in an algebraic closure of $F_p$.

3. $2^r \leq [K : F_p] \leq 2^{r+u-1}$

4. $\sigma_n \in \text{Gal}(K/F_p)$. $G_n$ is a cyclic subgroup of $K^\star$.

   Suppose $n = p^l d$.

5. $\sigma_d \in \text{Gal}(K/F_p)$. $G_n \subseteq G_{p^i d^j}$ for all $i, j \geq 0$.

6. There are integers $i_1, j_1, i_2, j_2$ such that $0 \leq i_1, j_1, i_2, j_2 \leq \sqrt{2^{n+u-1}}$, $(i_1, j_1) \neq (i_2, j_2)$ and such that $\sigma_{p^{i_1} d^{j_1}} = \sigma_{p^{i_2} d^{j_2}}$

7. $|G_n| / p^{i_1} d^{j_1} - p^{i_2} d^{j_2}$.

8. From the fact $2^{2^s} > n^{\sqrt{2^s}}$ it is easily deduced that for all $v \geq 0$, $2^{2^{s+v}} > n^{\sqrt{2^{s+2v}}}$. In particular, $2^{2^r} > n^{\sqrt{2^{r+u-1}}}$.

9. From the fact $(1 + m\theta) \in G_n$ for all $m \in \mathcal{S}$ we deduce as in Section 3 that $G_n$ contains $2^r$ different monomials over $F_p$. Hence, $|G_n| \geq 2^{2^r}$.

10. From items 6, 7, 8, 9 can be deduced that $d = 1$ so $n = p^l$.

11. Since $n$ is not a non trivial perfect power $n = p$. $\square$

**Corollary 5.1** If $n, k, a, u$ are as in the previous theorem then the primality of $n$ can be determined in $2^{2(u-1)} 2^{\max(s+2(u-1)-k,0)} \tilde{O}((\log n)^4)$ time.

**Proof**

As in the analysis of complexity of the previous sections. □

To be more precise about this result let $A_u$ the algorithm associated to Theorem 5.1 and $C(A_u)$ its complexity. Corollary 5.1 implies that $C(A_u) \approx 2^{4(u-1)} C(A_1)$ if $k \leq 2^s$ and $C(A_u) \approx 2^{2(u-1)} C(A_1)$ if $k \geq 2^{s+2(u-1)}$.

Even more precise, $C(A_u) \approx 2^{4(u-1)} 2^{-\min(\max(k-s,0),2(u-1))} C(A_1)$.

Note also that in the *unconditioned* case $(u = k)$ the complexity is $2^{4(k-1)} \tilde{O}((\log n)^6)$ which is polynomial time only for values of $k$ not too large.

## 5.2 The case $n \equiv -1 \pmod{4}$

Similarly when $n \equiv -1 \pmod{4}$ we have the following theorems, that we state without proof since the details are very similar to the previous results.

**Theorem 5.2** *Let $n \equiv -1 \pmod{4}$. Let $k = \nu_2(n+1)$. Let $s = \lceil 2 \log \log n \rceil$ and $t = s + 1$. Let $\alpha \in Z[i]$ and $u$ a positive integer, $1 \leq u \leq k+1$ and such that $\alpha^{\frac{n^2-1}{2^u}} \equiv -1 \pmod{n}$. Suppose that every positive integer less or equal than $2^{\max(s-k+2(u-1),0)+1}$ is coprime with $n$.*

*Suppose also that for every $m \leq 2^{\max(s-k+2(u-1),0)}$ we have $(1+mx)^n \equiv (1+mx^n) \bmod (n, x^{2^{t+2u-1}} - \alpha)$ and that $n$ is not a perfect power. Then, $n$ is prime.*

**Corollary 5.2** *If $n, k, \alpha, u$ are as in the previous theorem then the primality of $n$ can be determined in $2^{2u} 2^{\max(s+2(u-1)-k,0)} \tilde{O}((\log n)^4)$ time. In other words, if we call these tests $B_u$, then $C(B_u) \approx 4C(A_u)$.*

**Remark 5.1** *In each of the Theorems 3.1, 3.2, 4.1, 4.2, 5.1, 5.2, and Corollaries 5.1, 5.2, we claim that $s$ can be replaced by $\lceil 2 \log \log n \rceil$ (and $t$ by $\lceil 2 \log \log n \rceil + 1$). In fact, $s$ can be replaced by the minimum positive integer $s$ such that $|G_n| > n^{2^{s/2}}$. That $s \leq \lceil 2 \log \log n \rceil$ was achieved using the fact that $G_n$ contains properly the set $\mathcal{T}_1$ whose cardinality is $2^{2^s} - 1$. But actually $G_n$ contains the larger set $\mathcal{T}$ whose cardinality is the combinatorial number*

$$\binom{2^{s+1} - 1}{2^s} = \frac{1}{2} \binom{2^{s+1}}{2^s}.$$

*Using Stirling's formula with error, see for instance [9], it is easy to prove that $|G_n| > 2^{2^{s+1} - \frac{s}{2} - 3}$.*

*The smallest integer for which $2^{2^{s+1} - \frac{s}{2} - 3} > n^{2^{s/2}}$ is the smallest integer for which $2^{\frac{s}{2}+1} > \log n + s/2 + 3$. Since we know $s/2 \leq \lceil \log \log n \rceil$ then $s$ is at most the smallest value for which*

$$s/2 + 1 > \log(\log n + \lceil \log \log n \rceil + 3). \tag{1}$$

*It follows that $s = \lceil 2 \log \log n \rceil$ or maybe even $\lceil 2 \log \log n \rceil - 1$.*

*The algorithm should start by verifying which of the values satisfies (1) since each reduction in the value of $s$ in one unit improves around four times the speed of the algorithm.*

# 6 Conclusions and Conjecture

In practice, it is clearly desirable to apply algorithm 1 of Section 3 or algorithm 2 of Section 4 when possible. In the worst case $\nu_2(n-1) = k = 2$), algorithm 1 runs at least $2^{11}$ times faster than the best possible running time of the AKS algorithm for primes $n$ large enough. Hence, the worst case of algorithm 2 runs $2^9$ times faster than the best possible case of AKS. This occurs because the main step of Algorithm 1 executes at most $2^{s-2} \leq \frac{(\log n)^2}{4}$ iterations, each of which consist in multiplying polynomials of degree at most $(\log n)^2$. In contrast, in the best possible case AKS executes $8(\log n)^2$ multiplications of polynomials of degree at least $64(\log n)^2$. When $k$ is large the difference in the performance improves dramatically.

For implementation, if no integer $a$ satisfying $(\frac{a}{n}) = -1$ is known a priori, then a search for such an $a$ within a reasonable range should be implemented. In addition, if this fails to produce such an $a$, then a search for a small value of $u$ would be useful.

It is to be remarked that when the value of $k$ is small, the running time for these tests is still large. This indicates that it may be reasonable to develop analogous tests for numbers $n$ with large $\nu_m(n^f - 1)$ for reasonably small $f$.

Note that if $k > \frac{1}{2} \log n$ then the algorithms 1 and 2 run in $\tilde{O}(\log n)^2$ time. Also, while $k$ increases from 2 to $[2 \log \log n]$ the running time improves up to $\tilde{O}(\log n)^4$. But when $k$ varies from $[2 \log \log n]$ to $[\frac{1}{2} \log n]$ there is no more improvement in the speed of our algorithm. Here we believe one should attempt to sharpen the algorithms because the order of the group $G_n$ can be proven to increase together with $k$, in such a way that it forces $s$, the smallest solution of $|G_n| > n^{2^{s/2}}$, to decrease. To be precise we formulate the following conjecture, which we hope to prove in the near future.

**Conjecture.** *Algorithm 1 and 2 can be modified in such a way that while $k$ increases from 2 to $(1/2) \log n$ the complexity of both algorithms decreases from $\tilde{O}(\log n)^6$ to $\tilde{O}(\log n)^2$.*

**Acknowledgements.** *I am grateful to María González Lima for her help and support during the elaboration of this manuscript. I am grateful to Luis Báez Duarte for improving greatly the clarity of exposition of this paper and to Thomas G. Berry, for reading earlier versions and pointing out many errors. Finally, I am grateful to Víctor Ramírez for pointing out that $\sigma_d$ is also an automorphism and for other helpful discussions. His remark allowed to improve in 4 times the speed of the algorithms presented.*